\documentclass[preprint,12pt]{elsarticle}

\usepackage{amssymb}
\usepackage{amscd,amsmath,amsfonts,graphicx}
\usepackage{amsthm}

\newcommand{\E}{{\sf E}}
\newcommand{\e}{\mbox{e}}

\begin{document}
\baselineskip=20pt


\begin{center}
{\Large \bf Characterization of probability distributions via
 functional equations of power-mixture type}
\end{center}

\begin{center}
\vspace{0.4cm} Chin-Yuan Hu\footnote{National Changhua
University of Education, Taiwan, ROC} \quad 
Gwo Dong Lin\footnote{Hwa-Kang Xing-Ye Foundation and Academia Sinica,
Taiwan, ROC. E-mail: gdlin@stat.sinica.edu.tw} \quad
Jordan M. Stoyanov\footnote{Bulgarian Academy of Sciences, Bulgaria. E-mail: stoyanovj@gmail.com}
\end{center}

{\small
 \vspace{0.1cm}
\begin{center}
 {\bf Abstract} 
\end{center} 
\noindent
We study power-mixture type {\it functional
equations} in terms of Laplace--Stieltjes transforms of probability
distributions. These equations arise when studying {\it distributional
equations} of the  type \ $Z=  X + TZ$,
where $T$ is a known random variable, while the variable $Z$ is
defined via $X$, and we want to `find' $X$. We provide necessary and
sufficient conditions for such functional  equations to have  {\it
unique solutions}. The uniqueness is equivalent to a {\it characterization
property} of a probability distribution. We present results which are either new or
 extend and improve
 previous results about functional equations of compound-exponential and
compound-Poisson types. In particular, we give another affirmative answer to a
question posed by J. Pitman and M. Yor in 2003. We provide explicit
illustrative examples and deal with related topics. 
{\bf Running title:} Characterization of distributions via functional equations.\\
{\bf Corresponding author:} Jordan M. Stoyanov}

\vspace{0.4cm}\noindent{\bf 1. Introduction}

We deal with probability distributions and their characterization
properties expressed in the form of distributional
equations of the type  $Z  \stackrel{d}{=}  X + TZ$, where $T$ is a given random variable, the variable $Z$ is
defined via $X$, and we want to `find' $X$.
By using Laplace--Stieltjes transform (for short: LS-transform)
of the distributions of the random variables involved, we transfer such a distributional equation to a functional equation
of a specific type. Our goal is to provide necessary and
sufficient conditions for such a functional  equation to have a {\it
unique solution}. The unique solution is equivalent to a {\it characterization
property} of a probability distribution.

It is worth mentioning that the topic {\it Distributional Equations}
was intensively studied over the last decades. There are excellent sources; among them are 
the recent books by Buraczewski, Damek and Mikosch \cite{Bur16} and Iksanov \cite{Iks16}.
For good reasons, the phrase
 {\it ``The  equation} $X \ = \ AX + B"$ is included as a subtitle of \cite{Bur16}. 
From different perspectives this distributional equation is studied also in \cite{Iks16}.
Such equations are called `fixed-point equations'; they arise as
limits when studying autoregressive sequences in economics and
actuarial modelling, and the `fixed point' (the unique solution) is
related to the so-called perpetuities. These books contain a detailed
analysis of diverse stochastic models, a variety of results and
methods. Besides the authors of the two books, an essential contribution in this area 
is made by many scientists, to list here only a few names: H.
Kesten, C. Goldie, W. Vervaat, P. Embrechts, Z. Jurek, G. Alsmeyer. Much more can be found in the books cited above.

In the present paper, we study a wide class of power-mixture
functional equations for the LS-transforms of probability
distributions. In particular, equations of  compound-exponential
type, compound-Poisson type, and others, fall into this class. On
the other hand, the related   Poincar\'e type functional equations
have been studied by Liu \cite{Liu02} and recently by Hu and Lin
\cite{Hu20}; see also the references therein.

The power-mixture functional equations arise, e.g.,  when studying power-mixture transforms involving 
two sii-processes. Here the abbreviation `sii-processes' stands for a stationary-independent-increments stochastic processes.
Think, in particular, of L\'evy processes. Consider a continuous time sii-process $(X_1(t))_{t\ge 0},$
and let $F_{1,t}$ be the (marginal) distribution of $X_1(t);$ we write this as $X_1(t)\sim F_{1,t}.$
Moreover, let $X_1:=X_1(1)\ge 0$ be
the generating random variable for the process, so $X_1 \sim F_1:=F_{1,1}$  uniquely determines the distribution
 of the process $(X_1(t))_{t\ge 0}$ at any time $t$. Thus  we have the multiplicative semigroup
$(\hat{F}_{1,t}(s))_{t\ge 0}$ satisfying the power relation
\begin{eqnarray}
\hat{F}_{1,t}(s)=(\hat{F}_{1}(s))^t,\ \ s, t\ge 0.
\end{eqnarray}
Here
$\hat{F}_{1,t}$ is the LS-transform  of the
distribution $F_{1,t}$ of $X_1(t):$
\begin{eqnarray*}
\hat{F}_{1,t}(s)={\E}[{\e}^{-sX_1(t)}]=\int_0^{\infty}{\e}^{-sx}\,{\rm d}F_{1,t}(x),\  \ s \ge 0
\end{eqnarray*}
(see, e.g., Steutel and van Harn \cite{Ste04}, Chapter I).

Let further, $(X_2(t))_{t\ge 0},$  independent of $(X_1(t))_{t\ge 0},$ be another continuous time sii-process with a
generating random variable $X_2:=X_2(1)\ge 0$  and let $X_2(t)\sim
F_{2,t},\ X_2\sim F_2 := F_{2,1}.$ Now, we can consider the composition process
$(X(t))_{t\ge 0}:=(X_1(X_2(t))_{t\ge 0},$
which is the subordination of the process   $(X_1(t))_{t\ge 0}$ to the process $(X_2(t))_{t\ge 0}.$
The generating random variable for $(X(t))_{t\ge 0}$ is
$X:= X(1)=X_1(X_2(1))\sim F.$  In view of
Eq.\,(1), the distribution $F$ has LS-transform $\hat{F}$, which is of the power-mixture type (for short, power-mixture transform), and satisfies the following relations:
\begin{eqnarray}\hat{F}(s)&:=&{\E}[{\e}^{-sX}]=\int_0^{\infty}{\E}[{\e}^{-sX_1(u)}]\,{\rm d}F_{2}(u)
=\int_0^{\infty}(\hat{F}_{1}(s))^u\,{\rm d}F_{2}(u)\\
&=&\int_0^{\infty}\exp(-u\,[-\log
\hat{F}_1(s)])\,{\rm d}F_2(u)=\hat{F_2}(-\log \hat{F}_1(s)),\ \ s\ge 0.
\end{eqnarray}

 From now on, we will focus mainly on the power-mixture
transforms (2) or (3). The brief illustration of dealing with two
sii-processes is just one of the motivations. Thus, we now require
only the random variable $X_1 \sim F_1$ to be {\it infinitely
divisible}, but not asking this property for $X_2 \sim F_2$. For
such distributions $F$ with elegant LS-transforms, see Steutel and
van Harn \cite{Ste04}, Chapter III, as well as Pitman and Yor
\cite{Pit03}.

If $X_2 \sim F_2$, where $F_2 \in {\rm Exp}(1),$ the standard exponential distribution, \  $F_2(x)=1-{\e}^{-x},\
x \ge 0,$ its LS-transform is $\hat{F}_2(s)=1/(1+s), \ s\ge 0,$ and the distribution $F$ for the composition process
$(X(t))_{t \ge 0}$  reduces to the so-called {\it compound-exponential distribution} whose LS-transform
(for short, {\it compound-exponential transform}) is:
\begin{eqnarray}
\hat{F}(s)=\frac{1}{1-\log \hat{F}_1(s)},\ \ s\ge 0.
\end{eqnarray}
This shows that the power-mixture transforms are essentially more
general than the compound-exponential ones. The latter case, however, is
important by itself and it has been studied by Hwang and Hu \cite{Hwa11}.

When the random variable $X_1 \sim F_1$ is actually related to (or
constructed from) the variable $X \sim F,$ the LS-transform
$\hat{F}_1$ will be a function of the LS-transform $\hat{F}.$ Hence
the distribution $F$ (equivalently, its LS-transform $\hat{F}$) can
be considered as the solution to some of the functional equations
(2), (3) or (4). Since each of these equations is related to
 a distributional equation, as soon as we have a unique solution (a `fixed point'), this will
provide a characterization property of the corresponding distribution.

The main purpose of this paper is to provide necessary and
sufficient conditions for the functional equations in question to
have  {\it unique} distributional solutions, and we do this under general requirements.  
We exhibit new results; some of them either extend or  improve
previous results for functional equations of the
compound-exponential and compound-Poisson types. In particular, we provide another affirmative
answer to a question posed by Pitman and Yor \cite{Pit03}. This question and the answer were first 
given by Iksanov \cite{Iks02}, \cite{Iks04}. Our arguments are different; details are given in Example 2 below. 
Functional equations of other types are also studied.

In Section 2, we formulate the problem and state the main results
and corollaries. The results are  illustrated in  Section 3
by examples which fit well to the problem. Section 4 contains a
series of lemmas which we need in Section 5 for proving the main
theorems.  We conclude in Section 6 with comments and challenging
questions. The list of references includes significant works all related to our study.

\vspace{0.4cm}\noindent
{\bf 2. Formulation of the problem. Main results}

Let $X$ be a nonnegative random variable with distribution $F$ and
mean  $\mu={\E}[X],$  a number in the open interval $(0,\infty).$
Starting with $X \sim F,$ we will construct an {\it infinitely
divisible} random variable $X_1 \sim F_1$ to be used in Eq.\,(2).
Consider three nonnegative random variables and their distributions
as follows: \ $T \sim F_T, \ A \sim F_A, \ B \sim F_B.$ Suppose
further that $Z$ is a random variable, independent of $T,$ with the
{\it length-biased distribution} $F_Z$ induced by $F,$ namely,
\begin{eqnarray}
F_Z(z)=\frac{1}{\mu}\int_0^zx\,{\rm d}F(x),\ \ z\ge 0.
\end{eqnarray}
We involve also the scale-mixture random variable $TZ \sim F_{TZ}$. We are now prepared
to define the following two functions in terms of LS-transforms:
\begin{eqnarray}
\sigma(s)&:=&\mu\int_0^{s}\hat{F}_{TZ}(x)\,{\rm d}x=\int_0^{\infty}\frac{1-\hat{F}(ts)}{t}\,{\rm d}F_T(t),\ \ s\ge 0,\\
\sigma_B(s)&:=&\int_0^s\hat{F}_B(t)\,{\rm d}t,\ \ s\ge 0.
\end{eqnarray}
Notice that $\sigma ( \cdot )$ and $\sigma_B ( \cdot )$ are
Bernstein functions and their first derivatives  are completely
monotone functions, by definition; see, e.g., Schilling et
al.\,\cite{Sch12}. The function $\sigma$ in (6) will play a crucial
role in this paper and the integrand $(1-\hat{F}(ts))/t$ is defined
for $t=0$ by continuity to be equal to $\mu\,s.$ The second equality
in (6) can be verified by differentiating its both sides with
respect to $s$ and  using the following facts:
\[
\hat{F}_{Z}(s)=\E[\e^{-sZ}]=\frac{-\hat{F}'(s)}{\mu},\quad \int_0^s\hat{F}_Z(x){\rm d}x=\frac{1-\hat{F}(s)}{\mu}, \ \ s\ge 0.
\]
\indent Recall that in general the composition of two Bernstein
functions is a Bernstein function, hence this is so for
$\sigma_B\circ \sigma$, the functions in (6) and (7). We need also
the `simple' function, $\rho(s)=\e^{-s}, \ s \ge 0,$ which is the
LS-transform of the degenerate random variable at the point 1, and
use its property of being completely monotone. Therefore we can
consider the infinitely divisible random variable $X_1 \sim F_1$ (in
Eq.\,(1)) with LS-transform of compound-Poisson type: 
\begin{eqnarray}
\hat{F}_1(s)=\rho((\sigma_B\circ \sigma)(s))=\exp(-\sigma_B(\sigma(s))),\ \ s\ge 0.
\end{eqnarray}
Such a choice is appropriate in view of Lemmas 1 and 2 in Section 3. Clearly, $\hat{F}_1$ is a
function of $F, \ F_T$ and $F_B.$ Let us formulate our main results and some corollaries.

\vspace{0.2cm}\noindent
{\bf Theorem 1.}
{\it Under the above setting,  we have the following relations for $T, \ A$ and $B$:
\begin{eqnarray}
0 \le {\E}[T]<1, \quad {\E}[A]=1, \quad {\E}[A^2]<\infty\ \  {\rm and}  \quad 0 \le {\E}[B]<\infty,
\end{eqnarray}
if and only if the functional equation of power-mixture type
\begin{eqnarray}
\hat{F}(s)=\int_0^{\infty}\{\exp(-\sigma_B(\sigma(s)))\}^a\,{\rm d}F_A(a),\ \ s\ge 0,
\end{eqnarray}
has exactly one solution $X \sim F$ with mean $\mu$ and  finite variance. Moreover,
\begin{eqnarray}
{\sf Var}[X]=\frac{{\sf Var}[A]  +  {\E}[B]  +  {\E}[T]}{1  -  {\E}[T]}\,\mu^2.
\end{eqnarray}
}
\indent If we impose a condition on $B$, and use a.s. for `almost surely', Theorem 1 reduces as follows.

\vspace{0.1cm}\noindent {\bf Corollary 1.}
{\it In addition to the above setting, let $B=0 \ a.s.$ Then we have
\begin{eqnarray*}
0 \le {\E}[T]<1,\quad {\E}[A]=1\ \  {\rm and}  \quad {\E}[A^2]<\infty,
\end{eqnarray*}
if and only if the functional equation of power-mixture type
\begin{eqnarray*}
\hat{F}(s) =\int_0^{\infty}\exp\left(-a\,\sigma(s)\right)\,{\rm d}F_A(a),\ \ s\ge 0,
\end{eqnarray*}
has exactly one solution $X \sim F$ with mean $\mu$ and  finite variance. Moreover,
\begin{eqnarray*}
{\sf Var}[X]=\frac{{\sf Var}[A]  +  {\E}[T]}{1  -  {\E}[T]}\,\mu^2.
\end{eqnarray*}
}
\indent  If we impose a condition also on $A,$  Corollary 1 further reduces to the following.

\vspace{0.1cm}\noindent {\bf Corollary 2.}
{\it In addition to the setting in Theorem 1, let $A=1$\, a.s. and \ $B=0 \ a.s.$ Then
\[
0 \le {\E}[T]<1
\]
if and only if the functional equation of compound-Poisson type
\begin{eqnarray*}
\hat{F}(s)=\exp\left(-\sigma(s)\right),\ \ s\ge 0,
\end{eqnarray*}
has exactly one solution $X \sim F$ with mean $\mu$ and  finite variance. Moreover,
\begin{eqnarray*}
{\sf Var}[X]=\frac{{\E}[T]}{1  -  {\E}[T]}\,\mu^2.
\end{eqnarray*}
}
\indent
Here is a case of a `nice' proper random variable $A$, namely  $A \sim {\rm Exp}(1)$, so
$F_A(x)=1-\e^{-x}, \ x \ge 0.$  Corollary 1 takes now the following form.

\vspace{0.1cm}\noindent {\bf Corollary 3.}
{\it Let $X \sim F$
have mean  $\mu$, \ $B = 0$ a.s., \ $A \sim {\rm Exp}(1)$ and $T$ be a nonnegative random variable. Then
\[
0 \le{\E}[T]<1
\]
if and only if the  functional equation of compound-exponential type
\begin{eqnarray}
\hat{F}(s)=\frac{1}{1+ \sigma(s)}, \ s\ge 0,
\end{eqnarray}
has exactly one solution  $X \sim F$ with mean  $\mu$ and  finite variance. Moreover,
\begin{eqnarray}
{\sf Var}[X]=\frac{1  +  {\E}[T]}{1  -  {\E}[T]}\,\mu^2.
\end{eqnarray}
}
\indent And here is another particular but interesting case.

\vspace{0.1cm}\noindent {\bf Corollary 4.}
{\it In addition to the setting in Theorem 1, suppose that $T=p \ a.s.$ for some fixed  number $p \in (0,1)$  and that $B=0 \  a.s.$
Then we have
\[
{\E}[A]=1\ \quad {\rm and}  \quad {\E}[A^2]<\infty,
\]
if and only if the functional equation
\begin{eqnarray*}
\hat{F}(s)
=\int_0^{\infty}\exp\bigg(-a\,\frac{1-\hat{F}(ps)}{p}\bigg)\,{\rm d}F_A(a),\ \ s\ge 0,
\end{eqnarray*}
has exactly one solution $X\sim F$ with mean  $\mu$ and  finite variance. Moreover,
\begin{eqnarray*}
{\sf Var}[X]=\frac{{\sf Var}[A]  +  p}{1-p}\,\mu^2.
\end{eqnarray*}
} \indent  We now return to the construction of the infinitely
divisible LS-transform $\hat{F}_1$ in Eq. (8).  Using the completely
monotone function $\rho(s)=1/(1+\lambda s), \ s \ge 0$ (which
corresponds to ${\rm Exp}(\lambda)$), we have instead the LS-transform
\[
\hat{F}_1(s)=\rho((\sigma_B\circ \sigma)(s))=\frac{1}{1+\lambda\sigma_B(\sigma(s))},\ \ s \ge 0,
\]
and here is the next result.

\vspace{0.2cm}\noindent{\bf Theorem 2.} {\it Suppose, as before,
that $X \sim F$ is a nonnegative random variable with mean $\mu,$ a
number in the interval $(0,\infty)$. Let further \ $T, \ A$ and $B$
be three nonnegative random variables. Then, for a fixed constant $\lambda>0$, we have
\begin{eqnarray}
0 \le {\E}[T]<1, \quad {\E}[A]=1/\lambda,  \quad {\E}[A^2]<\infty\ \quad {\rm and} \quad 0 \le {\E}[B]<\infty,
\end{eqnarray}
if and only if the functional equation of power-mixture type
\begin{eqnarray}
\hat{F}(s)=\int_0^{\infty}\frac{1}{(1+\lambda\sigma_B(\sigma(s)))^{a}}\,{\rm d}F_A(a),\ \ s\ge 0,
\end{eqnarray}
has exactly one solution \ $X\sim F$ with mean $\mu$ and  finite variance. Moreover,
\begin{eqnarray}
{\sf Var}[X]=\frac{\lambda^2\,{\sf Var}[A]  +  \lambda+\E[B]  +  {\E}[T]}{1  -  {\E}[T]}\,\mu^2.
\end{eqnarray}
}
\indent Exchanging the roles of the arguments $a$ and $\lambda$ in Theorem 2 leads to  the following.

\vspace{0.2cm}\noindent{\bf Theorem 3.} {\it Consider the nonnegative random variables  $X, \ T, \ B, \ \Lambda$,
where $X \sim F$ has mean  $\mu,$ a positive number. Then, for an arbitrary constant $a>0,$ we have
\begin{eqnarray}
0 \le {\E}[T]<1, \quad {\E}[\Lambda]=1/a,\  \ {\E}[\Lambda^2]<\infty\ \ {\rm and}\ \ 0\le{\E}[B]<\infty,
\end{eqnarray}
if and only if the functional equation
\begin{eqnarray}
\hat{F}(s)=\int_0^{\infty}\left(1+\lambda
\sigma_B(\sigma(s))\right)^{-a}\,{\rm d}F_{\Lambda}(\lambda),\ \ s\ge 0,
\end{eqnarray}
has exactly one solution \ $X\sim F$ with mean  $\mu$ and  finite variance. Moreover,
\begin{eqnarray}
{\sf Var}[X]=\frac{a^2\,{\sf Var}[\Lambda]+a\,{\E}[\Lambda^2]+\E[B]+{\E}[T]}{1-{\E}[T]}\,\mu^2.
\end{eqnarray}
} \indent  In Theorems 2 and 3,  keeping both  $A$ and $\Lambda$ to
be proper random variables, that is, not a.s. constants,
 allows us to arrive at the following general result. For simplicity, $A$ and $\Lambda$ below are assumed to be independent.

\vspace{0.2cm}\noindent{\bf Theorem 4.} {\it Let   $X, \ T, \ A,\
\Lambda$ and $B$ be nonnegative random variables, where $X \sim F$
has mean  $\mu \in (0,\infty).$  We also require  $A$ and $\Lambda$ to be independent.  Then we have
\begin{eqnarray}
0 \le {\E}[T]<1, \quad {\E}[A\Lambda]=1, \quad {\E}[A^2]<\infty, \quad {\E}[\Lambda^2]<\infty \quad {\rm and} \quad 0\le {\E}[B]<\infty,
\end{eqnarray}
 if and only if the functional equation
\begin{eqnarray}
\hat{F}(s)=\int_0^{\infty}\int_0^{\infty}\left(1+\lambda\sigma_B(\sigma(s))\right)^{-a}\,{\rm d}F_{A}(a)\,{\rm d}F_{\Lambda}(\lambda), \ s\ge 0,
\end{eqnarray}
has exactly one solution $X\sim F$ with mean $\mu$ and  finite variance. Moreover,
\begin{eqnarray}
{\sf Var}[X]=\frac{{\sf Var}[A\Lambda]+{\E}[A\Lambda^2]+\E[B]+{\E}[T]}{1-{\E}[T]}\,\mu^2.
\end{eqnarray}
}

Clearly, when $\Lambda=\lambda=const \ a.s., $ Eqs.\,(20)--(22)
reduce to Eqs.\,(14)--(16), respectively, while  if $A=a=const \
a.s.,$ Eqs.\,(20)--(22) reduce to Eqs.\,(17)--(19), accordingly.
This is why in Section 5 we omit the proofs of Theorems 2 and 3,
however we provide a detailed proof of the more general Theorem 4.

Finally, let us involve the Riemann-zeta function defined as usual by
\[
\zeta(s)=\sum_{n=1}^{\infty}\frac{1}{n^s}, \ s>1.
\]
It is well known  that for any $a>1,$ the function
$\rho(s):={\zeta(a+s)}/{\zeta(a)}, \ s\ge 0,$ is the LS-transform of
a probability distribution which is called Riemann-zeta
distribution, and remarkably, it is infinitely divisible (see Lin
and Hu \cite{Lin01}, Corollary 1).  We have the following result
which is in the spirit of the previous theorems, however it is interesting by itself.

\vspace{0.2cm}\noindent{\bf Theorem 5.} {\it  Let $X, T$ and
$\Lambda$ be nonnegative random variables and $X \sim F$ have  mean
 $\mu$, a number in the interval $(0,\infty).$ Then, for any fixed
number $a>1$, we have
\begin{eqnarray}
0 \le {\E}[T]<1, \quad {\E}[\Lambda]=\frac{-\zeta(a)}{\zeta'(a)} \quad {\rm and} \quad {\E}[\Lambda^2]<\infty,
\end{eqnarray}
if and only if the functional equation
\begin{eqnarray}\hat{F}(s)=\frac{1}{\zeta(a)}\int_0^{\infty}\zeta(a+\lambda
\sigma(s))\,{\rm d}F_{\Lambda}(\lambda),\ \ s \ge 0,
\end{eqnarray}
has exactly one solution $X\sim F$ with mean  $\mu$ and  finite
variance. Moreover,
\begin{eqnarray}
{\sf Var}[X]=\frac{\zeta''(a)\,{\E}[\Lambda^2]-\zeta(a)+\zeta(a)\,{\E}[T]}{\zeta(a)(1-{\E}[T])}\,\mu^2.
\end{eqnarray}
}

\vspace{0.3cm} \noindent{\bf 3. Examples}

We present now some  examples to illustrate the use of the above
results. The first two examples can be considered as improvements of
Theorems 1.1 and 1.3 in Hwang and Hu \cite{Hwa11}. We use below the
notation  \ $\stackrel{{\rm d}}{=}$ \  meaning equality in distribution.

\vspace{0.1cm}\noindent {\bf Example 1.} We start with a random
variable $X$, where $0 \le X \sim F$ has  mean  $\mu \in
(0,\infty),$ and let $T$ be a nonnegative random variable. Assume
that $Z\ge 0$ is a random variable  with  the length-biased
distribution (5) induced by $F,$ and that $X_1, \ X_2$ are two
random variables each having the distribution $F.$
 Assume further that all random variables $Z, \ T, \ X_1, \ X_2$ are independent. Then
\[
0 \le {\E}[T]<1
\]
if and only if the distributional equation
\begin{eqnarray}
Z \ \stackrel{{\rm d}}{=} \ X_1+ X_2+T\,Z
\end{eqnarray}
has exactly one solution $X \sim F$ with mean $\mu$ and finite variance of the form (13).

This is true because the distributional equation (26) is equivalent
to the functional equation (12) expressed in terms of the
LS-transform $\hat{F}.$ Let us give details. We rewrite Eq.\,(26) as follows:
\[
\hat{F}_Z(s)=(\hat{F}(s))^2\,\frac{\sigma'(s)}{\mu},\ \ s\ge 0.
\]
By using the identity $\hat{F}_Z(s)=-\hat{F}'(s)/\mu,$ the above relation is equivalent to
\[
\frac{{\rm d}}{{\rm d}s}(\hat{F}(s))^{-1}=\sigma'(s),\ \ s\ge 0.
\]
 This means that indeed Eq.\,(12) holds true in view of the facts that $\hat{F}(0)=1$ and $\sigma(0)=0.$

Let us discuss two specific choices of \ $T,$  each one arriving at
interesting conclusion.

(a) When $T=0 \ a.s.,$ we have, by definition, $\sigma(s)=\mu s, \ s
\ge 0,$ and hence, by (12), $\hat{F}(s)=1/(1+\sigma(s))=1/(1+\mu s),
\ s \ge 0.$ Equivalently, $F$ is an exponential distribution with
mean  $\mu.$ On the other hand, Eq.\,(26) reduces to $Z
\stackrel{{\rm d}}{=} X_1+ X_2.$ Therefore, this equation claims to
be a characterization of the exponential distribution. The explicit
formulation is:

{\it The convolution of  an underlying distribution $F$ with itself is equal to the length-biased
distribution induced by $F$, if and only if, $F$ is an exponential distribution.}

(b) More generally, if $T=p \ a.s.$ for some fixed number $p \in
[0,1),$ then the unique solution $X \sim F$ to  Eq.\,(26) is the
following explicit mixture distribution
\[
F(x)=p+(1-p)(1 - \e^{-\beta x}), \ x \ge 0, \mbox{ where } \ \beta=(1-p)/\mu.
\]

\vspace{0.1cm}\noindent {\bf Example 2.} \ As in Example 1, we
consider two nonnegative random variables, $T$ and $X$, where $X
\sim F$ has mean  $\mu \in (0,\infty).$ Assume that the random
variable $Z \ge 0$ has  the length-biased distribution (5) induced
by $F,$ and that all random variables $X, \ T, \ Z$ are independent. Then
\[
0 \le {\E}[T]<1
\]
if and only if the distributional equation
\begin{eqnarray}
Z \ \stackrel{{\rm d}}{=} \ X+T\,Z
\end{eqnarray}
has exactly one solution $X \sim F$ with mean  $\mu$ and finite variance. Moreover,
\begin{eqnarray*}
{\sf Var}[X]=\frac{{\E}[T]}{1-{\E}[T]}\,\mu^2.
\end{eqnarray*}

Let us underline that this answers one of the questions posed by
Pitman and Yor \cite{Pit03}, p.\,320. The question itself can be
read (in our format) as follows:

{\it Given a random variable $T \sim F_T$, does there exist a random
variable $X \sim F$ (with unknown $F$) such that Eq.\,(27) is
satisfied with $Z$ having a length-biased distribution induced by
$F$? }

In order to explain the affirmative answer, note that the distributional equation (27) is equivalent to the functional
equation (by following the same idea as in Example 1):
\begin{eqnarray}
\hat{F}(s)={\e}^{-\sigma(s)}, \ s \ge 0.
\end{eqnarray}
This, however, is exactly the case of Corollary 2 (or, of Theorem 1
with $A=1 \ a.s.$ and $B=0 \ a.s.$).

It is seen that given any distribution of $T \ge 0$ with $0 \le
\E[T]<1,$ Eq.\,(27) characterizes the corresponding underlying
distribution (unique solution) $F$ of $X$ with mean $\mu$ and {\it
finite variance}. The behavior of the solution $F$ heavily depends
on the conditions on $T.$

Note that Iksanov \cite{Iks02, Iks16, Iks04} was the first to
provide an affirmative answer to the question by Pitman and Yor. His
conditions and conclusions are different from ours (the proofs are
of course different). For example, assuming that $T>0,$ $\E[\log T]$
exists (finite or infinite) and $\mu\in(0,\infty),$ Iksanov
\cite{Iks02} proved that there exists a unique solution $F$ (to
Eq.\,(27)) with mean $\mu$ if and only if $\E[\log T]<0;$ there is
no conclusion/condition about the variance of $F.$ Moreover, in our
condition ($0\le\E[T]<1$), we do not exclude the possibility that
${\sf P}[T=0]>0.$ Actually, it can be shown that if $T>0$ and
$\E[T]\in(0,1),$ then $\E[\log T]<0$ (because the function
$g(t)=t-1-\log t\ge 0$ for $t>0$). So if $T>0,$ our condition and
conclusion are stronger than those of Iksanov.

Let us consider four cases of $T.$

(a) If $T=0 \ a.s.,$ Eq.\,(27) reduces to $Z \stackrel{{\rm d}}{=}
X.$ It tells that  the length-biased distribution $F_Z$ is equal to
the underlying distribution $F.$  This equation characterizes the
degenerate distribution concentrated at the point $\mu$ because
Eq.\,(28) accordingly reduces to $\hat{F}(s)={\e}^{-\mu s}, \ s \ge
0.$

(b) If $T$ is a continuous random variable uniformly distributed on the interval $[0,1],$
 Eq.\,(27) characterizes the exponential distribution with mean  $\mu$ (see
also Pitman and Yor \cite{Pit03}, p.\,320). Indeed, by using the
identity
\[
\log(1+s)=\int_1^{\infty}\frac{s}{x(x+s)}\,{\rm d}x,\ \ s\ge 0,
\]
we see that the function $\hat{F}(s)=1/(1+\mu\,s), \ s\ge 0,$
satisfies Eq.\,(28).

More generally, if $T$ has a uniform distribution on the interval
$[p,1]$ for some $p \in [0,1),$ then the unique solution to
Eq.\,(28) is the following explicit mixture distribution
\[
F(x)=p+(1-p)(1-\e^{-\beta x}), \ x \ge 0, \ \mbox{ where } \ \beta=(1-p)/\mu.
\]

(c) If we assume now that $T$ has a beta distribution
$F_T(x)=1-(1-x)^{a}, \ x \in (0,1),$ with parameter $a>0,$ then the
unique solution $X \sim F$ to Eq.\,(27) will be the Gamma
distribution $F=F_{a,b}$ with density
\[
f_{a,b}(x)=\frac{1}{\Gamma(a)\,b^{a}}\,x^{a-1}\,{\e}^{-x/b},\ \ x>0.
\]
Here $b=\mu/a$ and we use the following identity:  \ for $a>0, \ b>0,$
\[
\log(1+ b\,s)=\int_0^1\frac{(1-t)^{a-1}}{t}[1-(1+ bst)^{-a}]\,{\rm d}t,\ \ s \ge 0,
\]
or, equivalently,
\[
\int_0^1\frac{a\,(1-t)^{a-1}}{(1+bst)^{a+1}}\,{\rm d}t=\frac{1}{1+b\,s},\ \  s \ge 0
\]
(see, e.g., Gradshteyn and Ryzhik \cite{Gra14}, Formula 8.380(7),
p.\,917).

 (d) Take a particular value $\mu=2/3$ and assume that $T$ has the density
$g(t)=1/\sqrt{t}-1,\ t \in (0,1).$ Then Eq.\,(27) has a unique
solution $X \sim F$ with LS-transform
$\hat{F}(s)=2s/(\sinh\sqrt{2s})^2, \ s> 0$ (expressed in terms of
the hyperbolic-sine function; see Pitman and Yor \cite{Pit03},
p.\,318). In general, if $\mu \in (0,\infty)$ is an arbitrary number
(not specified) and $T$ is as above, then the unique solution $X
\sim F$ has LS-transform $\hat{F}(s)=3\mu s/(\sinh\sqrt{3\mu s})^2,
\ s >0.$

Notice that Eq.\,(27) can also be solved by fitting to the
Poincar\'e type functional equation considered in Theorem 4 of Hu
and Lin \cite{Hu20}. This idea, however, requires the third moment
of the underlying distribution $F$ to be involved.

On the other hand,  we can replace $Z$ in Eq.\,(27) by a random
variable $X^*$ which obeys the {\it equilibrium distribution} $F^*$
induced by $F.$ Recall that
\begin{eqnarray}
F^*(x)=\frac{1}{\mu}\int_0^x\bar{F}(t)\,{\rm d}t,\ \ x\ge 0,
\end{eqnarray}
where $\bar{F}(t)={\sf P}[X>t]=1-F(t), \ t \ge 0.$ In this case we obtain an
interesting characterization result, and this is the content of the next example.

\vspace{0.1cm}\noindent {\bf Example 3.}  Let $0 \le X\sim F$ with
mean  $\mu \in (0,\infty)$ and let $T$ be a nonnegative random
variable. Assume that the random variable  $X^* \sim F^*$ follows
the equilibrium distribution defined in (29). Further, assume that
all random variables $X, \ T, \ X^*$ are independent. Then we have
\[
0 \le{\E}[T]<1
\]
if and only if the distributional equation
\begin{eqnarray}
X^* \ \stackrel{{\rm d}}{=} \ X+T\,X^*
\end{eqnarray}
has exactly one solution $X\sim F$ with mean  $\mu$ and finite
variance of the form (13).

Indeed, this is true because the distributional equation (30) is
equivalent to the functional equation (12). The latter follows from
rewriting Eq.\,(30) in terms of LS-transforms:
\begin{eqnarray}
\hat{F}_{X^*}(s)&=&\hat{F}(s)\,{\E}[{\e}^{-sTX^*}]
=\hat{F}(s)\int_0^{\infty}{\E}[{\e}^{-stX^*}]\,{\rm d}F_T(t)\nonumber\\
&=&\hat{F}(s)\int_0^{\infty}\hat{F}_{X^*}(st)\,{\rm d}F_T(t),\ \ s \ge 0.
\end{eqnarray}
We need to use also the relation \
$\hat{F}_{X^*}(s)=(1-\hat{F}(s))/(\mu\,s), \ s> 0$ (see Lemma 8(ii)
below). Plugging this identity in (31) and carrying out the function
$\hat{F}$ leads to Eq.\,(12).

As before, letting \ $T=0 \ a.s.$ in (30), we get another
characterization of the exponential distribution (because, by (12),
$\hat{F}(s)=1/(1+\mu\,s), \ s \ge 0$). The statement is:

{\it The equilibrium distribution $F^*$ (see (29) above) is equal to
the underlying distribution $F,$  if and only if, $F$ is
exponential.} (See also Cox \cite{Cox62}, p.\,63.)

\vspace{0.2cm}\noindent{\bf 4. Ten Lemmas}

To prove the main results, we need some auxiliary statements given
here as lemmas. The first two lemmas are well known and Lemma 1 is
called {\it Bernstein's Theorem} (see, e.g., Steutel and van Harn
\cite{Ste04}, p.\,484, or Schilling et al.\,\cite{Sch12}, p.\,28).

\vspace{0.2cm}\noindent{\bf Lemma 1.}
{\it The LS-transform $\hat{F}$ of a nonnegative random variable $X\sim F$ is a completely
monotone function on $[0,\infty)$ with $\hat{F}(0)=1,$ and vice versa.}

\vspace{0.2cm}\noindent{\bf Lemma 2.}
{\it (a) \ The class of Bernstein functions  is closed under composition. Or, the composition of
two Bernstein functions is still  a Bernstein function.\\
(b) \ Let \ $\rho$ \ be a completely monotone function and  $\sigma$ a Bernstein function on $[0,\infty).$ Then their composition
$\rho\circ\sigma$ is a completely monotone function on $[0,\infty).$}

 Note that in Theorems 1 and 2 we have used two simple choices for the function $\rho.$
 The next two lemmas concern the contraction property of some `usual' real-valued functions of real arguments.
These properties will be used later to prove the uniqueness of the solution to functional equations in question.

\vspace{0.2cm}\noindent{\bf Lemma 3.}
{\it  Let $a, b\ge 0.$ Then:\\
(i) \ $|\log(1+a)-\log(1+b)|\le |a-b|;$\\
(ii) \ $|{\rm e}^{-a}-{\rm e}^{-b}| \le |a-b|.$}

\vspace{0.1cm}\noindent
{\bf Proof.} Since $a$ and $b$ are exchangeable, it is enough to show the validity of (i) and (ii)
for  $a \ge b \ge 0.$ 
For claim (i), consider the function $g(x)=\log (1+x)- x,\ x \ge 0.$ Since
$g'(x)=(1+x)^{-1}-1 \le 0$ for  $x \ge 0,$ $g$ is a decreasing
function on $[0,1].$ Therefore, $g(a) \le g(b)$ for $a\ge b\ge 0.$ Equivalently, $\log(1+a)- \log(1+b)\le a-b,$ and hence,
\[
|\log(1+a)-\log(1+b)|=\log(1+a)-\log(1+b)\le a-b=|a-b|, \ a \ge b\ge 0.
\]
For claim (ii), we use the inequality
${\e}^{-b}-{\e}^{-a}=\int_b^a {\e}^{-t}\,{\rm d}t \le \int_b^a1\,{\rm d}t=a-b, \ a \ge b\ge 0.$
Therefore, $|{\e}^{-a}-{\e}^{-b}|={\e}^{-b}-{\e}^{-a} \le a-b=|a-b|, \ a \ge b \ge 0.$ The proof is complete.

\vspace{0.2cm}\noindent{\bf Lemma 4.}
{\it  (i) \ For arbitrary \ $a, \ b \in [0,1]$ \ and \ $t\ge 1,$  we have:
\[
|a^t-b^t| \le t\,|a-b|.
\]
(ii) \ For real numbers \ $x, y \ge 0$ and $a>1,$  the Riemann-zeta function satisfies
\[
|\zeta(a+x)-\zeta(a+y)|\le -\zeta'(a)\,|x-y|.
\]
(iii) \ For any $a>1,$ we have \ $\zeta''(a)\zeta(a)>(\zeta'(a)^2.$
}

\vspace{0.1cm}\noindent {\bf Proof.} It is easy to establish claim
(i); still, details can be seen in Hu and Lin \cite{Hu20}. For claim
(ii), we use Lemma 3(ii). Indeed,
\begin{eqnarray*}
&~&|\zeta(a+x)-\zeta(a+y)|=\left|\sum_{n=1}^{\infty}\frac{1}{n^{a+x}}-\sum_{n=1}^{\infty}
\frac{1}{n^{a+y}}\right|\\
&\le &\sum_{n=1}^{\infty}\frac{1}{n^{a}}\left|\frac{1}{n^x}-\frac{1}{n^y}\right|=\sum_{n=1}^{\infty}\frac{1}{n^{a}}\left|{\rm e}^{-x\log n}-{\rm e}^{-y\log n}\right|\\
&\le &\sum_{n=1}^{\infty}\frac{1}{n^{a}}\left|x\log n-y\log n\right|=\sum_{n=1}^{\infty}\frac{\log n}{n^a}|x-y|=-\zeta'(a)\,|x-y|.
\end{eqnarray*}
We have used the fact that $\zeta'(s)=-\sum_{n=1}^{\infty}(\log
n)/n^s$ for $s>1.$ To prove claim (iii), we consider the nonnegative
random variable $X$ whose LS-transform is
\[\pi(s)=\zeta(a+s)/\zeta(a), \ s\ge 0.\] Then $\E[X]=-\lim_{s\to
0^+}\pi'(s)=-\zeta'(a)/\zeta(a)$ and $\E[X^2]=\lim_{s\to
0^+}\pi''(s)=\zeta''(a)/\zeta(a)$ (see Lemma 6 below). The required
inequality follows from the fact that ${\sf Var}[X] ={\E}[X^2] -
({\E}[X])^2>0.$ The proof is complete.

\vspace{0.2cm}
We need now  notations for the first two moments of the random variable $X \sim F$ and a useful relation
implied by the positivity of the variance ${\sf Var}[X]$:
\[
m_1={\E}[X],   \quad m_2={\E}[X^2] \quad \mbox{ with } \quad  m_1^2 \leq m_2.
\]
Notice that instead of `first moment $m_1$', sometimes it is
convenient to use the equivalent `mean  $\mu$', as we have already
done.

\vspace{0.2cm}\noindent{\bf Lemma 5.} {\it  Suppose the nonnegative
random variable $X\sim F$ has finite positive second moment.  Then
the LS-transform $\hat{F}$ has a sharp upper bound as follows:
\begin{eqnarray}
\hat{F}(s) \le 1-\frac{m_1^2}{m_2}+\frac{m_1^2}{m_2}\,{\rm e}^{-(m_2/m_1)s},\ \ s\ge 0.
\end{eqnarray}}
\indent For the proof of Lemma 5 we refer to Eckberg \cite{Eck77},
Guljas et al.\,\cite{Gul98} or Hu and Lin \cite{Hu08}. It is
interesting to mention that the RHS of the inequality (32) is
actually the LS-transform of a specific two-point random variable
$X_0 \sim F_0$ (with first two moments $m_1,m_2$). Indeed, define
the values of $X_0$ and their probabilities
 as follows: \[{\sf P}[X_0=0]=1-\frac{m_1^2}{m_2}\quad \ {\rm and}~\quad \ {\sf P}[X_0=\frac{m_2}{m_1}]=\frac{m_1^2}{m_2}.\]

\vspace{0.1cm} Here is another result, Lemma 6; its proof is given
in Lin \cite{Lin93}.

\vspace{0.2cm}\noindent{\bf Lemma 6.} {\it  Let $0 \le X\sim F$ with
LS-transform $\hat{F}.$ Then for each integer $n \ge 1,$ the $n$th
order moment of $X$, finite or infinite, can be calculated as
follows:
\[
m_n := {\E}[X^n]=\lim_{s \to 0^+}(-1)^n\hat{F}^{(n)}(s)=(-1)^n\hat{F}^{(n)}(0^+).
\]
} \indent Let us deal again with equilibrium distributions. For a
random variable $X$, \ $0 \le X \sim F$ with finite positive mean
$\mu$ (= first moment $m_1$),
 we define the first-order equilibrium distribution based on $F$ by
$F_{(1)}(x) := \frac{1}{\mu}\int_0^x{\bar F}(y)\,{\rm d}y, \ x \ge
0$ (in Eq.\,(29), we have used the notation $F^*$).  If we assume
that for some $n,$ \ $m_n = {\E}[X^n] < \infty,$ we define
iteratively the equilibrium distribution $F_{(k)}$ of order $k,$ for
\ $k= 1, 2, \ldots, n,$ as follows:  $F_{(k)}(x) :=
\frac{1}{\mu_{(k-1)}}\int_0^x{\bar F}_{(k-1)}(y)\,{\rm d}y, \ x \ge
0.$ We have used here the notation $\mu_{(j)}$ for the mean (the
first moment) of $F_{(j)}$: \ $\mu_{(j)} := \int_0^{\infty} x\,{\rm
d}F_{(j)}(x)$. Also, with \ $F_{(0)}={F}$, \ $\mu_{(0)}= m _1,$ \
$m_0=1$ (the total mass is 1), we achieve full consistency.

It is clear from the above definition that $m_n < \infty$ implies
that $\mu_{(n-1)} < \infty,$ and vice versa. Moreover, finite are
all moments $m_k$ and all means $\mu_{(k)}$ for $k < n.$

We state in Lemma 7 below an interesting relationship between the
means  $\{\mu_{(k)}\}$ and the moments $\{m_k\}$. For details see,
e.g., Lin \cite{Lin98}, p.\,265, or Harkness and Shantaram
\cite{Har69}.

\vspace{0.2cm}\noindent{\bf Lemma 7.} {\it  Let for some integer $n
\ge 2$ the $n$th order moment $m_n$ of the random variable $0 \le
X\sim F$ be strictly positive and finite.   Then, for any $k=1, 2,
\ldots, n-1$, the mean  $\mu_{(k)}$ of the $k$th-order equilibrium
distribution $F_{(k)}$ is well defined (finite) and moreover,
\[
\mu_{(k-1)}=\frac{m_k}{k\,m_{k-1}} \quad \mbox{ for } \quad k=1, 2,  \ldots, n.
\]
} \indent For the proofs of the last three lemmas, we refer to  Hu
and Lin \cite{Hu20}.

\vspace{0.2cm}\noindent{\bf Lemma 8.} {\it  Consider the nonnegative
random variable $X \sim F$ whose mean $\mu$ is strictly positive and
finite,
 and let $X^{*} \sim F^{*},$ where $F^*$ is the equilibrium distribution induced by $F.$
 Then for $s>0,$ the following statements are true:\\
(i) \ $(1-\hat{F}(s))/s=\int_0^{\infty}{\rm e}^{-sx}(1-F(x))\,{\rm d}x;$\\
(ii) \ $\hat{F}^{*}(s)=(1-\hat{F}(s))/(\mu s)\le 1;$\\
(iii) \ $(\hat{F}(s)-1+\mu\,s)/s^2=\mu\int_0^{\infty}{\rm e}^{-sx}(1-F^{*}(x))\,{\rm d}x;$\\
(iv) \ $\lim_{s\to 0^+}(1-\hat{F}(s))/s=\mu;$\\
(v) \ $\lim_{s \to 0^+}(\hat{F}(s)-1+\mu s)/s^2=\frac12\,{\E}[X^2]$
\ (finite or infinite).}

\vspace{0.2cm}\noindent{\bf Lemma 9.} {\it  Given is a sequence of
random variables \ $\{Y_n\}_{n=1}^{\infty}$, where $Y_n \geq 0$ and
$Y_n \sim G_n.$ We impose two assumptions:\\ (a) all \ $Y_n$, hence
all $G_n$, have the same finite first two moments, that is, \
${\E}[Y_n]=m_1, \ {\E}[Y_n^2]=m_2$ \ for \  $n=1, 2, \ldots;$ \\ (b)
the LS-transforms $\{\hat{G}_n\}_{n=1}^{\infty}$ form a decreasing
sequence of functions.\\ Then the following limit exists:
\[
\lim_{n\to\infty}\hat{G}_n(s) =: \hat{G}_{\infty}(s), \ s\ge 0.
\]
Moreover, \ $\hat{G}_{\infty}$ is the LS-transform of the
distribution $G_{\infty}$ of a random variable $Y_{\infty} \geq 0$
with first moment  $\E[Y_{\infty}]=m_1$ and second moment
$\E[Y_{\infty}^2]$ belonging to the interval \ $[m_1^2, m_2].$}

\vspace{0.2cm}\noindent{\bf Lemma 10.} {\it  Suppose that \ $W_1
\sim F_{W_1}$ and \ $W_2 \sim F_{W_2}$ are nonnegative random
variables with the same mean  (same first moment) $\mu_W$, a
strictly positive finite number. Consider another random variable
$Z_* \geq 0$, \ where $Z_* \sim F_{Z_*}$ has a positive mean
$\mu_{Z_*}<1.$ Assume further that the LS-transforms of \ $W_1$ and
$W_2$ satisfy the following relation:
\begin{eqnarray}
|\hat{F}_{W_1}(s)-\hat{F}_{W_2}(s)|\le\int_0^{\infty}|\hat{F}_{W_1}(ts)-\hat{F}_{W_2}(ts)|\,{\rm d}F_{Z_*}(t), \ s\ge 0,
\end{eqnarray} or, equivalently,
\[
\big|\E[{\rm e}^{-sW_1}]-\E[{\rm e}^{-sW_2}]\big| \le
\big|\E[{\rm e}^{-sZ_*W_1}]-\E[{\rm e}^{-sZ_*W_2}]\big|,\  \ s\ge 0.
\]
Then \ $\hat{F}_{W_1}=\hat{F}_{W_2}$ and hence $F_{W_1}=F_{W_2}.$}

\vspace{0.4cm}\noindent{\bf 5. Proofs of the main results}

We start with the proof of Theorem 1, then omit details about
Theorems 2 and 3, however provide the proof of the more general
Theorem 4. Finally we give the proof of Theorem 5. Each of the
proofs consists naturally of two steps, Step 1 (Sufficiency) and
Step 2 (Necessity). In many places, in order to make a clear
distinction between factors in long expressions, we use the dot
symbol, `` $\cdot$ ", for multiplication.

\vspace{0.2cm}\noindent{\bf Proof of Theorem 1.}

\vspace{0.1cm}\noindent {\it Step 1 (Sufficiency).}  Suppose that
Eq.\,(10) has exactly one solution, $X$, where $0 \le X \sim F$ with
mean  $\E[X] = \mu \in(0,\infty)$ and  finite variance (and hence
$\E[X^2] < \infty).$ Then we want to prove that all conditions (9)
are satisfied.


First, rewrite Eq.\,(10) as follows:
\[
\hat{F}(s)=\int_0^{\infty}\exp\bigg(-a\int_0^{\sigma(s)} \hat{F}_B(t)\,{\rm d}t\bigg)\,{\rm d}F_A(a),\ \ s \ge 0.
\]
Differentiating twice this relation with respect to $s,$ we find, for $s>0,$ that
\begin{eqnarray}
\hat{F}'(s)&=&\int_0^{\infty}(-a)\exp\bigg(-a\int_0^{\sigma(s)} \hat{F}_B(t)dt\bigg)\,{\rm d}F_A(a)\cdot\hat{F}_B(\sigma(s))\sigma'(s),\\
\hat{F}''(s)&=&\int_0^{\infty}a^2\exp\bigg(-a\int_0^{\sigma(s)} \hat{F}_B(t)\,{\rm d}t\bigg)\,{\rm d}F_A(a)\cdot(\hat{F}_B(\sigma(s))\sigma'(s))^2\nonumber\\
& &+\int_0^{\infty}(-a)\exp\bigg(-a\int_0^{\sigma(s)} \hat{F}_B(t)\,{\rm d}t\bigg)\,{\rm d}F_A(a)\cdot\hat{F}^{\prime}_B(\sigma(s))(\sigma'(s))^2\nonumber\\
& &+\int_0^{\infty}(-a)\exp\bigg(-a\int_0^{\sigma(s)} \hat{F}_B(t)\,{\rm d}t\bigg)\,{\rm d}F_A(a)\cdot\hat{F}_B(\sigma(s))\sigma''(s).
\end{eqnarray}
Letting \ $s\to 0^+$ \ in (34) and (35) yields, respectively,
\begin{eqnarray*}
\hat{F}'(0^+)&=&\hat{F}'(0^+)\,\E[A],\\
\hat{F}'(0^+)&=&\E[A^2](\hat{F}'(0^+))^2-\E[A]\left(\hat{F}'_B(0^+)(\hat{F}^{\prime}(0^+))^2-\hat{F}''(0^+)\E[T]\right).
\end{eqnarray*}
Equivalently, in view of Lemma 6, we obtain two relations:
\begin{eqnarray}
\mu&=&\mu\,{\E}[A],\\
\E[X^2]&=&\E[A^2]\,\mu^2+{\E[A]}\,(\E[B]\,\mu^2+{\E}[X^2]\,E[T]).
\end{eqnarray}
\indent Since $\mu$ and $\E[X^2]$ are strictly positive and finite, we conclude
from (36) and (37) that $\E[A]=1$  and that each of the quantities $\E[A^2], \ \E[B], \ \E[T]$ is finite.
Moreover, $\E[T] \le 1$ due to (37) again.  We need, however, the strong inequality $\E[T]<1.$
Suppose on the contrary, namely  that $\E[T]=1.$ Then this would imply that
$\E[A^2]=0$ by (37), a contradiction to the fact that $\E[A]=1.$ This
proves that the conditions (9) are satisfied. In addition, relation (11) for the variance ${\sf Var}[X]$ also follows
from (36) and (37) because
\[
\E[X^2]=\frac{\E[A^2]+\E[B]}{1-\E[T]}\,\mu^2.
\]
The sufficiency part is established.

\vspace{0.2cm}\noindent
{\it Step 2 (Necessity).} Suppose now that the conditions (9) are satisfied. Then we will show the existence of
 a solution $X\sim F$ to Eq.\,(10) with mean $\mu$ and  finite variance.

 To find such a solution $X\sim F,$ we first {\it define} two numbers:
\begin{eqnarray}
m_1=\mu \quad \mbox{and} \quad  m_2=\frac{\E[A^2]+\E[B]}{1-\E[T]}\,m_1^2,
\end{eqnarray}
and show later these happen to be the first two moments of the
solution.   Note that the denominator $1-\E[T]>0$ by (9) and that
the  numbers $m_1,\,m_2$ do satisfy the required moment relation
$m_2 \ge m_1^2,$ because
 $\E[A^2] \ge (\E[A])^2=1$ due to (9) and  Lyapunov's inequality.
 Therefore, the RHS of (32) with
$m_1, \ m_2$ as expressed in  (38) is a {\it bona fide}
LS-transform, say $\hat{F}_0,$ of a nonnegative random variable $Y_0
\sim F_0$ (by Lemma 1). Namely,
\[
\hat{F}_0(s)=1-\frac{m_1^2}{m_2}+\frac{m_1^2}{m_2}\,{\e}^{-(m_2/m_1)s},\ \ s\ge 0.
\]
It is clear that $m_1,\,m_2$ are exactly the first two moments of
$Y_0 \sim F_0,$ as mentioned before.

Next, using the initial $Y_0\sim F_0$ we define iteratively a
sequence of random variables $\{Y_n\}_{n=1}^{\infty},$ $Y_n \sim
F_n,$ through their LS-transforms (see Lemma 2):
\begin{eqnarray}
\hat{F}_n(s)=\int_0^{\infty}\exp\bigg(-a\int_0^{\sigma_{n-1}(s)} \hat{F}_B(t)\,{\rm d}t\bigg)\,{\rm d}F_A(a),\ \ s\ge 0, \ n\ge 1,
\end{eqnarray}
where
\[
\sigma_{n-1}(s)=\int_0^{\infty}\frac{1-\hat{F}_{n-1}(ts)}{t}\,{\rm d}F_T(t),\ \ s\ge 0.
\]
Differentiating (39) twice with respect to $s$ and letting \ $s \to
0^+,$ we have, for $n \ge 1,$
\begin{eqnarray}
\hat{F}'_n(0^+)&=&\hat{F}'_{n-1}(0^+)\,\E[A],\\
\hat{F}''_n(0^+)&=&\E[A^2]\,(\hat{F}'_{n-1}(0^+))^2-
\E[A]\left(\hat{F}'_B(0^+)(\hat{F}'_{n-1}(0^+))^2-\hat{F}''_{n-1}(0^+)\,\E[T]\right).
\end{eqnarray}

By Lemma 6, induction on $n$ and  in view of  (40) and (41), we can
show that for any $n= 1, 2, \ldots$, we have $\E[Y_n]=\E[Y_0]=m_1$
and $\E[Y_n^2]=\E[Y_0^2]=m_2$ (see relations (38)). Hence,
\begin{eqnarray}
{\sf Var}[Y_n]=
m_2-m_1^2=\frac{{\sf Var}[A]+{\E}[B]+{\E}[T]}{1-{\E}[T]}\,m_1^2,\ \ n \ge 1.
\end{eqnarray}
Moreover, by Lemma 5, we first have $\hat{F}_1 \le \hat{F}_0$, and
then by the iteration (39), $\hat{F}_n\le \hat{F}_{n-1}$ for any $n
\ge 2.$  Namely, $\{Y_n\}_{n=0}^{\infty}$ is a sequence of
nonnegative random variables having the same first two moments $m_1,
m_2,$ and the LS-transforms $\{\hat{F}_n\}$ are decreasing.
Therefore, Lemma 9 applies. Denote the limit of $\{\hat{F}_n\}$ by
$\hat{F}_{\infty}.$ Then $\hat{F}_{\infty}$  will be the
LS-transform of the distribution, $F_{\infty}$, of a nonnegative
random variable, $Y_{\infty}$, that is, $Y_{\infty} \sim
F_{\infty}$, where $\E[Y_{\infty}]=m_1$ and $\E[Y_{\infty}^2] \in
[m_1^2, m_2].$ Thus it follows from (39) that the limit $F_{\infty}$
is a solution to Eq.\,(10) with mean $\mu=m_1$ and  finite variance.
Applying once again Lemma 6 to Eq.\,(10) (with $X=Y_{\infty}$ and
$F=F_{\infty}$), we conclude that $\E[Y_{\infty}^2]=m_2$ as
expressed in (38), and hence the solution ${Y}_{\infty} \sim
F_{\infty}$ has the required variance as shown in (11) or (42).

Finally, let us establish the uniqueness of the solution to Eq.\,(10). Suppose, under conditions (9), that there are two
solutions, say $X \sim F$ and $Y\sim G,$ each satisfying Eq.\,(10) and each having
 mean  equal $\mu$ (and hence both having the same finite variance
as shown above). Thus we want to show that $F=G,$ or, equivalently, that $\hat{F}=\hat{G}.$ Let us introduce two functions,
\[
\bar{\sigma}_F(s)=\int_0^{\infty}\frac{1-\hat{F}(ts)}{t}\,{\rm d}F_T(t), \quad
\bar{\sigma}_G(s)=\int_0^{\infty}\frac{1-\hat{G}(ts)}{t}\,{\rm d}F_T(t), \ s \ge 0.
\]
Then we have, by assumption,
\[
\hat{F}(s)=\int_0^{\infty}\exp\left(-a\sigma_B(\bar{\sigma}_F(s))\right)\,{\rm d}F_A(a), \quad
\hat{G}(s)=\int_0^{\infty}\exp\left(-a\sigma_B(\bar{\sigma}_G(s))\right)\,{\rm d}F_A(a),\ \ s \ge 0.
\]
 Using Lemma 3, we get the inequalities:
\begin{eqnarray*}
&~&|\hat{F}(s)-\hat{G}(s)|\le\int_0^{\infty}a\left|\int_0^{\bar{\sigma}_F(s)}
\hat{F}_B(t)\,{\rm d}t-\int_0^{\bar{\sigma}_G(s)}\hat{F}_B(t)\,{\rm d}t\right|\,{\rm d}F_A(a)\\
&\le&{\E}[A]\left|\int_{\bar{\sigma}_G(s)}^{\bar{\sigma}_F(s)}
\hat{F}_B(t)\,{\rm d}t\right|=\left|\int_{\bar{\sigma}_G(s)}^{\bar{\sigma}_F(s)}\hat{F}_B(t)\,{\rm d}t\right|
\le\left|\bar{\sigma}_F(s)-\bar{\sigma}_G(s)\right|,\ \ s\ge 0.
\end{eqnarray*}
We have used the fact that ${\E}[A]=1$.  Thus we obtain that for $s>0,$
\begin{eqnarray*}
\left|\frac{1-\hat{F}(s)}{\mu s}-\frac{1-\hat{G}(s)}{\mu s}\right|
&\le&\left|\int_0^{\infty}\frac{1-\hat{F}(ts)}{\mu ts}\,{\rm d}F_T(t)-\int_0^{\infty}\frac{1-\hat{G}(ts)}{\mu ts}\,{\rm d}F_T(t)\right|\\
&\le& \int_0^{\infty}\left|\frac{1-\hat{F}(ts)}{\mu ts}-\frac{1-\hat{G}(ts)}{\mu ts}\right|\,{\rm d}F_T(t).
\end{eqnarray*}
This relation is equivalent to another one, for the pair of
distributions $F^*$ and $G^*$, induced, respectively, by $F$ and
$G;$ see Lemma 8. Thus
\[
|\hat{F}^*(s)-\hat{G}^*(s)|\le\int_0^{\infty}\left|\hat{F}^*(ts)-\hat{G}^*(ts)\right|\,{\rm d}F_T(t),\ \ s>0.
\]
However this is exactly relation (33). Therefore, Lemma 10 applies,
because ${\E}[T]<1$ and $F^*$,  $G^*$ have the same mean  by Lemma
7. Hence $\hat{F}^*=\hat{G}^*$, which in turn implies that
$\hat{F}=\hat{G}$ due to the fact that $F$ and $G$ have the same
mean  (see Huang and Lin \cite{Hua95}, Proposition 1). The proof of
the necessity and hence of Theorem 1 is complete.

\vspace{0.2cm}\noindent{\bf Proof of Theorem 4.}
Although the proof has some similarity to that of Theorem 1, it is given here for completeness and reader's convenience.

\vspace{0.1cm}\noindent {\it Step 1 (Sufficiency).}  Suppose that
Eq.\,(21) has exactly one solution $0 \le X \sim F$ with mean $\mu,$
a positive and finite number,  and  finite variance (hence $\E[X^2]
\in (0,\infty)$). Now we want to show that all five conditions in
(20) are satisfied.

Differentiating twice Eq.\,(21) with respect to $s,$ we have, for $s>0,$ the following:
\begin{eqnarray}
\hat{F}'(s)&=&\int_0^{\infty}\int_0^{\infty}(-a)\lambda
\left(1+\lambda\sigma_B(\sigma(s))\right)^{-(a+1)}\,{\rm d}F_{A}(a)\,{\rm d}F_{\Lambda}(\lambda)\cdot\sigma'_B(\sigma(s))\sigma'(s),\\
\hat{F}''(s)&=&\int_0^{\infty}\int_0^{\infty}(-a)(-a-1)\lambda^2
\left(1+\lambda\sigma_B(\sigma(s))\right)^{-(a+2)}\,{\rm d}F_{A}(a)\,{\rm d}F_{\Lambda}(\lambda)\cdot(\sigma'_B(\sigma(s))\sigma'(s))^2\nonumber\\
&~&+\int_0^{\infty}\int_0^{\infty}(-a)\lambda
\left(1+\lambda\sigma_B(\sigma(s))\right)^{-(a+1)}\,{\rm d}F_{A}(a)\,{\rm d}F_{\Lambda}(\lambda)\cdot\sigma''_B(\sigma(s))(\sigma'(s))^2\nonumber\\
&~&+\int_0^{\infty}\int_0^{\infty}(-a)\lambda\left(1+\lambda\sigma_B(\sigma(s))\right)^{-(a+1)}\,{\rm d}F_{A}(a){\rm d}F_{\Lambda}(\lambda)\cdot
\sigma'_B(\sigma(s))\sigma''(s).
\end{eqnarray}
Letting \ $s \to 0^+$ in (43) and (44) yields, respectively,
\begin{eqnarray*}
\hat{F}'(0^+)&=&\hat{F}'(0^+)\E[A\Lambda],\\
\hat{F}''(0^+)&=&\E[A(A+1)\Lambda^2]\,(\hat{F}'(0^+))^2+\E[A\Lambda]\left(\E[B](\hat{F}'(0^+))^2+\E[T]\,\hat{F}''(0^+)\right).
\end{eqnarray*}
Equivalently, we have, by Lemma 6,
\begin{eqnarray}
\mu&=&\mu\ \E[A\Lambda],\\
\E[X^2]&=&\E[A(A+1)\Lambda^2]\,\mu^2+\E[A\Lambda]\left(\E[B]\,\mu^2+\E[X^2]\,\E[T]\right).
\end{eqnarray}
\indent From (45) and (46) it follows that  $\E[A\Lambda]=1$  and
that each of the quantities $\E[(A\Lambda)^2],$ \ $\E[\Lambda^2],$ \
$\E[B],$ \ $\E[T]$  is strictly positive and finite; this is because
$\mu$ and  $\E[X^2]$ are numbers in $(0,\infty).$ Moreover, $\E[T]
\le 1$ due to (46), and it remains to show the strict bound
$\E[T]<1.$ Suppose on the contrary that $\E[T]=1.$ Then we would
have $\E[(A\Lambda)^2]=0$ by (46), a contradiction to the fact that
$\E[A\Lambda]=1.$ Thus we conclude that all conditions in (20) are
satisfied. Besides, the expression for the variance ${\sf Var}[X]$
(see (22)) also follows from (45) and (46), because
\[
\E[X^2]=\frac{\E[A(A+1)\Lambda^2]+\E[B]}{1-\E[T]}\mu^2.
\]
The sufficiency part is established.

\vspace{0.1cm}\noindent
{\it Step 2 (Necessity).}  Suppose now that the conditions (20) are satisfied. We want to show the existence of
 a solution $X \sim F$ to Eq.\,(21) with mean $\mu$ and  finite variance.

 Set first
\begin{eqnarray}
m_1=\mu \quad \mbox{and} \quad m_2=\frac{\E[A(A+1)\Lambda^2]+\E[B]}{1-\E[T]}\,m_1^2.
\end{eqnarray}
As in the proof of Theorem 1, we have $1-\E[T]>0$, $m_2 \ge m_1^2$
and the existence of a nonnegative random variable (we use the same
notations) $Y_0 \sim F_0$, where the LS-transform $\hat{F}_0$ is
equal to the RHS of (32).
The next is to use the initial $Y_0\sim F_0$ and define iteratively the
sequence of random variables $Y_n \sim F_n,\ n=1,2,\ldots,$ through the LS-transforms (see Lemma 2):
\begin{eqnarray}
\hat{F}_n(s)=\int_0^{\infty}\int_0^{\infty}\left(1+\lambda\sigma_B(\sigma_{n-1}(s))\right)^{-a}\,
{\rm d}F_A(a)\,{\rm d}F_{\Lambda}(\lambda),\ \ s \ge 0, \ n \ge 1,
\end{eqnarray}
where
\[\sigma_{n-1}(s)=\int_0^{\infty}\frac{1-\hat{F}_{n-1}(ts)}{t}\,{\rm d}F_T(t),\ \ s\ge 0.
\]
Differentiating (48) twice with respect to $s$ and letting  $s \to
0^+,$ we have, for $n\ge 1,$
\begin{eqnarray}
\hat{F}'_n(0^+)&=&\hat{F}'_{n-1}(0^+)\E[A\Lambda],\\
\hat{F}''_n(0^+)&=&\E[A(A+1)\Lambda^2](\hat{F}'_{n-1}(0^+))^2\nonumber\\
&&+\ \E[A\Lambda]\left(\E[B]\,(\hat{F}'_{n-1}(0^+))^2+\E[T]\,\hat{F}''_{n-1}(0^+)\right).
\end{eqnarray}

By Lemma 6 and induction on $n,$ we find through (49) and (50) that
$\E[Y_n]=\E[Y_0]=m_1$ and $\E[Y_n^2]=\E[Y_0^2]=m_2$  for any $n \ge
1,$ and hence
\begin{eqnarray}
{\sf Var}[Y_n]=
m_2-m_1^2=\frac{{\sf Var}[A\Lambda]+{\E}[A\Lambda^2]+\E[B]+{\E}[T]}{1-{\E}[T]}\,m_1^2,\ \ n\ge 0.
\end{eqnarray}
Moreover, by Lemma 5, we first have $\hat{F}_1\le \hat{F}_0$, and
then by the iteration (48), $\hat{F}_n\le \hat{F}_{n-1}$ for all
$n\ge 2.$  Thus, $\{Y_n\}_{n=0}^{\infty}$ is a sequence of
nonnegative random variables having all the same first two moments
$m_1, m_2,$ such that the sequence of their LS-transforms
$\{\hat{F}_n\}$ is decreasing. Therefore, Lemma 9 applies, and the
limit  $\lim_{n \to \infty}\hat{F}_n =: \hat{F}_{\infty}$  is the
LS-transform of a nonnegative random variable $Y_{\infty} \sim
F_{\infty}$ with mean  $\E[Y_{\infty}]=m_1$ and second moment
$\E[Y_{\infty}^2] \in [m_1^2,m_2].$ Consequently, it follows from
(48) that the limit $F_{\infty}$ is a solution to Eq.\,(21) with
mean $\mu=m_1$ and  finite variance. Applying Lemma 6 to Eq.\,(21)
again (with $X=Y_{\infty}$ and $F=F_{\infty}$), we conclude that
$\E[Y_{\infty}^2]=m_2$ (as in (47)), and hence the solution
${Y}_{\infty} \sim F_{\infty}$ has the required variance as shown in
(22) or (51).

Finally, let us show the uniqueness of the solution to Eq.\,(21).
Suppose that, under conditions (20),  there are two solutions, $X
\sim F$ and $Y \sim G,$ which satisfy Eq.\,(21) and both have the
same mean  $\mu$ (hence the same finite variance).

Now we want to show that $F=G,$ or, equivalently, that
$\hat{F}=\hat{G}.$ We need the functions
\[
\bar{\sigma}_F(s)=\int_0^{\infty}\frac{1-\hat{F}(ts)}{t}\,{\rm d}F_T(t),
\quad \bar{\sigma}_G(s)=\int_0^{\infty}\frac{1-\hat{G}(ts)}{t}\,{\rm d}F_T(t), \ s\ge 0.
\]
Then we have
\[\hat{F}(s)=\int_0^{\infty}\int_0^{\infty}\left(1+\lambda\sigma_B(\overline{\sigma}_F(s))\right)^{-a}\,{\rm d}F_{A}(a)\,{\rm d}F_{\Lambda}(\lambda), \ s\ge 0,\]
\[\hat{G}(s)=\int_0^{\infty}\int_0^{\infty}\left(1+\lambda\sigma_B(\overline{\sigma}_G(s))\right)^{-a}\,{\rm d}F_{A}(a)\,{\rm d}F_{\Lambda}(\lambda), \ s\ge 0.\]

 Using Lemma 3, we obtain the following chain of the relations:
\begin{eqnarray*}
&~&|\hat{F}(s)-\hat{G}(s)| \le \int_0^{\infty}\int_0^{\infty}a\lambda\left|\int_0^{\bar{\sigma}_F(s)}\hat{F}_B(t)\,{\rm d}t-
\int_0^{\bar{\sigma}_G(s)}\hat{F}_B(t)\,{\rm d}t\right|\,{\rm d}F_A(a)\,{\rm d}F_{\Lambda}(\lambda)\\
&\le&{\E}[A\Lambda]\left|\int_{\bar{\sigma}_G(s)}^{\bar{\sigma}_F(s)}\hat{F}_B(t)\,{\rm d}t\right|
=\left|\int_{\bar{\sigma}_G(s)}^{\bar{\sigma}_F(s)}\hat{F}_B(t)\,{\rm d}t\right|
\le\left|\bar{\sigma}_F(s)-\bar{\sigma}_G(s)\right|, \ s\ge 0,
\end{eqnarray*}
where we have used the condition ${\E}[A\Lambda]=1.$  The remaining
arguments are similar to those in the proof of Theorem 1, so we can omit the details. Thus
the necessity is established and the proof of Theorem 4 is complete.

\vspace{0.2cm}\noindent{\bf Proof of Theorem 5.} We follow a similar
idea as in the proofs of Theorems 1 and 4. It will be convenient for
the reader to see the details, all explicitly expressed in terms of
the Riemann-zeta function.

\vspace{0.1cm}\noindent {\it Step 1 (Sufficiency).} Suppose that
Eq.\,(24) has exactly one solution $0 \le X \sim F$ with mean  $\mu
\in (0,\infty)$ and  finite variance (hence $\E[X^2] \in
(0,\infty)$). Thus we want to show that conditions (23) are
satisfied.

Differentiating twice Eq.\,(24) with respect to $s,$ we have, for $s>0,$
\begin{eqnarray}
\hat{F}'(s)&=&\frac{1}{\zeta(a)}\int_0^{\infty}\lambda\zeta'(a+\lambda
\sigma(s))\,{\rm d}F_{\Lambda}(\lambda)\cdot\sigma'(s),\\
\hat{F}''(s)&=&\frac{1}{\zeta(a)}\int_0^{\infty}\lambda^2\zeta''(a+\lambda
\sigma(s))\,{\rm d}F_{\Lambda}(\lambda)\cdot(\sigma'(s))^2\nonumber\\
&&+\ \frac{1}{\zeta(a)}\int_0^{\infty}\lambda\zeta'(a+\lambda
\sigma(s))\,{\rm d}F_{\Lambda}(\lambda)\cdot\sigma''(s).
\end{eqnarray}
Letting \ $s\to 0^+$ in (52) and (53) yields, respectively,
\begin{eqnarray*}
\hat{F}'(0^+)&=&\frac{-\zeta'(a)}{\zeta(a)}\hat{F}^{\prime}(0^+)\,\E[\Lambda],\\
\hat{F}''(0^+)&=&\frac{\zeta''(a)}{\zeta(a)}(\hat{F}^{\prime}(0^+))^2\,\E[\Lambda^2]-\frac{\zeta'(a)}{\zeta(a)}\,\hat{F}''(0^+)\E[\Lambda]\,\E[T].
\end{eqnarray*}
Equivalently, we have, by Lemma 6, the following relations:
\begin{eqnarray}
\mu&=&\mu\ \frac{-\zeta'(a)}{\zeta(a)}\,\E[\Lambda],\\
\E[X^2]&=&\frac{\zeta''(a)}{\zeta(a)}\,\E[\Lambda^2]\,\mu^2-\frac{\zeta'(a)}{\zeta(a)}\,\E[X^2]\,\E[\Lambda]\,\E[T].
\end{eqnarray}
\indent From (54) and (55) it follows that
$\E[\Lambda]={-\zeta(a)}/{\zeta'(a)}$  and that both quantities
$\E[\Lambda^2]$ and $\E[T]$ are finite, because $\mu, \ \E[X^2] \in
(0,\infty).$ Also, $\E[T] \le 1$ due to (55) again.  To prove the
strict inequality $\E[T]<1,$ we assume on the contrary  that
 $\E[T]=1.$ In such a case $\E[\Lambda^2]=0$ by (55), which contradicts the fact
$\E[\Lambda]={-\zeta(a)}/{\zeta'(a)}.$ Thus conditions (23) are satisfied. Besides, relation (25) also follows from (54) and (55) because
\[
\E[X^2]=\frac{{\zeta''(a)}\E[\Lambda^2]}{\zeta(a)(1-\E[T])}\,\mu^2.
\]
The sufficiency part is established.

\vspace{0.2cm}\noindent
{\it Step 2 (Necessity).}  Suppose that conditions (23) are satisfied. We want to show the existence of
 a solution $X \sim F$ to Eq.\,(24) with mean $\mu$ and finite variance.

 We start with the relations
\begin{eqnarray}
m_1=\mu \quad \mbox{and} \quad m_2=\frac{{\zeta''(a)}\,\E[\Lambda^2]}{\zeta(a)(1-\E[T])}\,m_1^2.
\end{eqnarray}
In (56) the denominator $1-\E[T]$ is strictly positive by (23) and
$m_2 \ge m_1^2,$ because
$\E[\Lambda^2]\ge(\E[\Lambda])^2=(\zeta(a)/\zeta'(a))^2\ge\zeta(a)/\zeta''(a)$
(see Lemma 4). Therefore, as before the RHS of (32) with $m_1, m_2$
as expressed
 in (56) is an LS-transform, say $\hat{F}_0,$ of a nonnegative random variable $Y_0 \sim F_0$ (by Lemma 1).
Thus, starting with $Y_0 \sim F_0$ we can define iteratively the
sequence of random variables $Y_n\sim F_n, \ n=1,2,\ldots,$ through
LS-transforms (see Lemma 2):
\begin{eqnarray}
\hat{F}_n(s)=\frac{1}{\zeta(a)}\int_0^{\infty}\zeta(a+\lambda
\sigma_{n-1}(s))\,{\rm d}F_{\Lambda}(\lambda), \ s\ge 0, \ n\ge 1,
\end{eqnarray}
where
\[
\sigma_{n-1}(s)=\int_0^{\infty}\frac{1-\hat{F}_{n-1}(ts)}{t}\,{\rm d}F_T(t),\ \ s\ge 0.
\]
Differentiating (57) twice with respect to $s$ and letting \ $s\to 0^+,$ we find, for $n\ge 1,$
\begin{eqnarray}
\hat{F}'_n(0^+)&=&\frac{-\zeta'(a)}{\zeta(a)}\,\hat{F}'_{n-1}(0^+)\,\E[\Lambda],\\
\hat{F}''_n(0^+)&=&\frac{\zeta''(a)}{\zeta(a)}\,(\hat{F}'_{n-1}(0^+))^2\,\E[\Lambda^2]-
\frac{\zeta'(a)}{\zeta(a)}\,\hat{F}''_{n-1}(0^+)\,\E[\Lambda]\E[T].
\end{eqnarray}

By Lemma 6, induction on $n$ and relations (58) and (59), we find that $\E[Y_n]=\E[Y_0]=m_1,$
$\E[Y_n^2]=\E[Y_0^2]=m_2$ (as expressed in (56)) for any $n \ge 1$ and hence
\begin{eqnarray}
{\sf Var}[Y_n]=
m_2-m_1^2=\frac{\zeta''(a)\,{\E}[\Lambda^2]-\zeta(a)+\zeta(a)\,{\E}[T]}{\zeta(a)\,(1-{\bf E}[T])}\,m_1^2, \ n \ge 0.
\end{eqnarray}
Moreover, by Lemma 5, we first have $\hat{F}_1\le \hat{F}_0$, and
then by the iteration (57), $\hat{F}_n\le \hat{F}_{n-1}$ for all $n\ge 2.$  Thus, $\{Y_n\}_{n=0}^{\infty}$ is a sequence of
nonnegative random variables having the same first two moments
$m_1, m_2,$ and decreasing sequence of their LS-transforms $\{\hat{F}_n\}$. Applying Lemma 9, there is a limit
 $\lim_{n \to \infty} \hat{F}_n =: \hat{F}_{\infty},$ which is the
LS-transform of a nonnegative random variable $Y_{\infty} \sim
F_{\infty}$ with $\E[Y_{\infty}]=m_1$ and $\E[Y_{\infty}^2] \in
[m_1^2, m_2].$ Hence, it follows from (57) that  $F_{\infty}$ is a
solution to Eq.\,(24) with mean $\mu=m_1$ and finite variance.
Applying again Lemma 6 to Eq.\,(24), with $X=Y_{\infty}$ and
$F=F_{\infty}$, we conclude that $\E[Y_{\infty}^2]=m_2$ as in (56),
and hence the solution ${Y}_{\infty} \sim F_{\infty}$ has the
required variance as shown in (25) or (60).

Finally, it remains to prove the uniqueness of the solution to
Eq.\,(24). Suppose, under conditions (23), that there are two
solutions, $X \sim F$ and $Y \sim G,$ satisfying Eq.\,(24) and
having the same mean  $\mu$ (hence the same finite variance). We
want to show that $F=G,$ or, equivalently, that $\hat{F}=\hat{G}.$
We use the functions
\[
\bar{\sigma}_F(s)=\int_0^{\infty}\frac{1-\hat{F}(ts)}{t}\,{\rm d}F_T(t), \quad
\bar{\sigma}_G(s)=\int_0^{\infty}\frac{1-\hat{G}(ts)}{t}\,{\rm d}F_T(t), \ s\ge 0,
\]
to express explicitly the two LS-transforms:
\[
\hat{F}(s)=\frac{1}{\zeta(a)}\int_0^{\infty}\zeta(a+\lambda
\bar{\sigma}_F(s))\,{\rm d}F_{\Lambda}(\lambda), \ s\ge 0,
\]
\[
\hat{G}(s)=\frac{1}{\zeta(a)}\int_0^{\infty}\zeta(a+\lambda
\bar{\sigma}_G(s))\,{\rm d}F_{\Lambda}(\lambda),\ \ s\ge 0.
\]

By Lemma 4, we derive the relations:
\begin{eqnarray*}
&~&|\hat{F}(s)-\hat{G}(s)|\le\frac{-\zeta'(a)}{\zeta(a)}\int_0^{\infty}\lambda\left|\bar{\sigma}_F(s)-\bar{\sigma}_G(s)\right|
\,{\rm d}F_{\Lambda}(\lambda)\\
&\le&\frac{-\zeta'(a)}{\zeta(a)}\,{\E}[\Lambda]\left|\bar{\sigma}_F(s)-\bar{\sigma}_G(s)\right|=\left|\overline{\sigma}_F(s)-\bar{\sigma}_G(s)\right|,\ \ s\ge 0,
\end{eqnarray*}
in which we have used the fact that
${\E}[\Lambda]={-\zeta(a)}/{\zeta'(a)}.$   The remaining arguments
are similar to those in the proofs of Theorems 1 and 4, and hence
omitted here. The necessity is established and  the proof of Theorem
5 is completed.

\vspace{0.2cm}\noindent{\bf 6. Concluding remarks}

Below are some useful remarks regarding the problems and the results
in this paper and their relations with previous works.

\vspace{0.2cm}\noindent{\bf Remark 1.} In Theorem 1, we have treated
the power-mixture type functional equation (Eq.\,(2)) which includes
the compound-Poisson equation, Eq.\,(28), as a special case. Thus,
the problems and the results here can be considered as an extension
of the previous works, in particular, the well known work by Pitman
and Yor \cite{Pit03}. We have given all necessary details in Example
2.

\vspace{0.1cm}\noindent{\bf Remark 2.} In Examples 1 and 3, when
$T=p\ a.s.$ for some fixed number \ $p\in[0,1),$ the unique solution
$X \sim F$ to Eqs.\,(26) and (30) is the mixture distribution
\[
F(x)=p+(1-p)(1-{\rm e}^{-\beta x}), \ x\ge 0, \ \mbox{ with } \ \beta=(1-p)/\mu.
\]
Its LS-transform has a mixture form:
\[
\hat{F}(s)=1-\frac{\mu}{\lambda}+\frac{\mu}{\lambda}\frac{1}{1+\lambda s}, \ s \ge 0, \
\mbox{ where } \ \lambda=\frac{1}{\beta}=\frac{\mu}{1-p}.
\]
 Actually,  for any $T \sim F_T$ with $0\le T\le 1,\ \E[T]<1$ and for any number
$p \in [0,1)$ such that  $F_T(p) \in (0,1],$ the unique solution $X \sim F$ to
Eqs.\,(26) and (30) satisfies the inequality:
\[
\hat{F}(s)\le 1-\frac{\mu}{\lambda}+\frac{\mu}{\lambda}\frac{1}{1+\lambda s}, \ s \ge 0, \
\mbox{ where } \ \lambda=\frac{\mu}{F_T(p)(1-p)}.
\]
Notice that this relation is satisfied even if the explicit form of
$\hat{F}$ is unknown.

\vspace{0.2cm}\noindent{\bf Remark 3.} The class of power-mixture
transforms defined in Eq.\,(2) is quite rich and includes the
LS-transforms of the so-called \ $C_t, \ S_t, \ T_t$ random
variables (where $t>0$), which
  are expressed in terms of the hyperbolic functions,
$\cosh, \sinh, \tanh,$ respectively.

Indeed, for \ $C_t, \ t>0$, we have the LS-transform:
\begin{eqnarray*}
\E[{\e}^{-sC_t}]&=&\bigg(\frac{1}{\cosh \sqrt{2s}}\bigg)^t=\exp(-t\log
(\cosh\sqrt{2s}))\\
&=&\exp\bigg(-t\int_0^s\frac{\tanh{\sqrt{2x}}}{\sqrt{2x}}dx\bigg),\
\ s>0.\end{eqnarray*} This is exactly the form of Eq.\,(2), where
 $X_2=t \ a.s.$ and $X_1 \sim F_1$ with
\[
\hat{F}_1(s)=\exp\bigg(-\int_0^s\frac{\tanh{\sqrt{2x}}}{\sqrt{2x}}dx\bigg),\ \ s\ge 0.
\]

 Similar arguments apply to the LS-transforms of the random variables
\ $S_t$ and  $T_t$:
\[
\E[{\e}^{-sS_t}]=\bigg(\frac{\sqrt{2s}}{\sinh \sqrt{2s}}\bigg)^t,\ \ s>0, \quad \mbox{and} \quad
\E[{\e}^{-sT_t}]=\bigg(\frac{\tanh \sqrt{2s}}{\sqrt{2s}}\bigg)^t,\ \ s> 0.
\]
 It is also interesting to note that for any fixed $t>0$, the following relation holds:
\[
\E[{\e}^{-sC_t}]=\E[{\e}^{-sS_t}]\,\E[{\e}^{-sT_t}],\ \ s> 0.
\]
Therefore, we have an interesting distributional equation
\[
C_t \ \stackrel{{\rm d}}{=} \ S_t+T_t.
\]
This means that the random variable $C_t$ can be decomposed into a
sum of two subindependent random variables  $S_t$ and $T_t.$
(See Pitman and Yor \cite{Pit03}.)

\vspace{0.1cm}\noindent{\bf Remark 4.} We finally consider an
equation which is similar to Eq.\,(27) (or Eq.\,(28)), but not
really of the power-mixture type  Eq.\,(10). Let $0\le X\sim F$ with
mean $\mu \in (0,\infty)$ and let $T$ be a nonnegative random
variable. Assume that the random variable $Z \ge 0$ has  the
length-biased distribution (5) induced by $F.$ Let the random
variables $X_1, X_2$ be independent copies of  $X \sim F,$ and
moreover, let $X_1,X_2,T$ be independent. Then the distributional
equation
\begin{eqnarray}
Z \ \stackrel{{\rm d}}{=} \ X_1+T\,X_2
\end{eqnarray}
(different from Eq.\,(27)) is equivalent to the functional equation
\[
\hat{F}(s)={\rm e}^{-\sigma_*(s)}, \ s\ge 0
\]
(compare with Eq.\,(28)).  Here the Bernstein function $\sigma_*$
 is of the form:
\[
\sigma_*(s)=\int_0^s\bigg(\mu\int_0^{\infty}\hat{F}(xt)\,{\rm d}F_T(t)\bigg)\,{\rm d}x,\ \ s\ge 0.
\]


To analyze the solutions to this kind of equations is a serious problem. The attempt to
follow the approach in this paper was not successful. Perhaps a new idea is needed.
However, there is a specific case when the solution to the above
equation is explicitly known. More precisely, let us take  $T \
\stackrel{{\rm d}}{=} \ U^2$ with $U$ being uniformly distributed on
$[0,1].$ In this case Eq.\,(61) has a unique solution: the
hyperbolic-cosine distribution, say $F,$ with LS-transform
\[
\hat{F}(s)=\bigg(\frac{1}{\cosh \sqrt{\mu\,s}}\bigg)^2,\ \ s\ge 0.
\]
Therefore \ $X \ \stackrel{{\rm d}}{=} \ \frac12\,{\mu}\,C_2$ (see,
e.g., Pitman and Yor \cite{Pit03}, p.\,317). Once again, this
characteristic property (Eq.\,(61) with $T \ \stackrel{{\rm d}}{=} \
U^2$) is found just for $C_t$ with $t=2$. What about arbitrary
$t>0$\,? As far as we know, for general random variables $C_t, \
S_t, \ T_t$ the characterizations of
 their distributions are challenging but still open problems.

\vspace{0.3cm}\noindent
{\bf Acknowledgments} \\
We are grateful to Sergey Foss, Gerold Alsmeyer and Alexander Iksanov for their attention and useful comments on previous 
version of our paper.

\end{document}